\documentclass[12pt]{article}
\usepackage{amssymb}
\usepackage{amsmath}
\newtheorem{theo}{Theorem}

\newtheorem{prop}[theo]{Proposition}
\newtheorem{lemm}[theo]{Lemma}
\newtheorem{coro}[theo]{Corollary}

\def \C{{\mathbb C}}
\def \R{{\mathbb R}}
\def \Z{{\mathbb Z}}

\def \p{{\mathfrak p}}
\def \Q{{\mathbb Q}}
\def \N{{\mathbb N}}
\def \F{{\mathbb F}}
\def \A{{\mathbb A}}

\begin{document}

\title{Existence of non-preperiodic algebraic points for a rational self-map
of infinite order}

\author{Ekaterina Amerik}

\date{June 18, 2010}

\maketitle

Let $X$ be a smooth projective variety defined over a number field $K$ 
 and let $f:X\dasharrow X$ be a dominant  
rational self-map defined over the same number field. As shown 
in \cite{AC}, one can attach to $f$ a dominant
rational map $g:X\dasharrow T$, commuting with $f$ and such that the fiber 
of $g$ through a sufficiently general complex point $x\in X(\C)$ is the
Zariski closure of its iterated orbit (or ``$f$-orbit'') 
$\{f^k(x),\ k\in \N\}$. Here ``sufficiently general'' means ``outside a
countable union of proper subvarieties'', and so this theorem does not give
any information on the $f$-orbits of algebraic points, which, apriori,
can have smaller Zariski closure than general complex points.

One would of course like to show that in reality this never happens
and one can always find an algebraic point such that its $f$-orbit
 is ``as large''
as the general one. For instance, a conjecture already  implicit in \cite{AC}
and formulated by Medvedev and Scanlon in \cite{MS} (Conjecture 5.3)
states that if no power of $f$ preserves a non-trivial fibration, then there 
should be a point $x\in X(\bar \Q)$ with Zariski-dense $f$-orbit ; a variant of this is an 
earlier 
conjecture by
S.-W. Zhang stating the same in the case when $f$ is regular and polarized
(that is, there is an ample line bundle $L$ on $X$ with $f^*L=qL$ for some
$q>1$).

What is certainly true in the case when $f$ is regular and polarized is
that, at least, there exist points in $X(\bar \Q)$ with infinite
$f$-orbits (that is, non-preperiodic algebraic points). The reason is
that in this case, one can introduce the so-called {\it canonical height}
$\hat{h}_L:X(\bar \Q)\to \R$  which is a Weil height function for $L$ 
with the property $\hat{h}_L(f(x))=q\hat{h}_L(x)$; it follows that
the set of preperiodic points is a set of bounded height and therefore
it cannot exhaust $X(\bar \Q)$ (see \cite{CS}). However, the theory of 
canonical heights
does not seem to work well enough for rational self-maps.

The purpose of this note is to provide an elementary proof of the 
existence of non-preperiodic algebraic points for such dominant rational
self-maps (using, though, a result by E. Hrushovski which does not
seem to have been treated in a very accessible way at the moment).
The argument is very similar to the one used by Bell, Ghioca and
Tucker to prove a version of the ``dynamical Mordell-Lang conjecture''
for unramified endomorphisms of quasiprojective varieties; in fact
this note is directly inspired by \cite{BGT} (and is, in some sense, a continuation of 
\cite{ABR} where some general observations about iteration of algebraic points
by a dominant rational self-map with a fixed point have been made). 
The point is 
that, thanks to Hrushovski's result, for some positive integer $k$ 
one can find an $f^k$-invariant
{\it $\p$-adic neighbourhood} in $X$, for a suitable prime $\p$
in the ring of integers of a suitable finite extension of $K$. Then
one uses \cite{BGT} to conclude that all preperiodic points in this 
neighbourhood are periodic with bounded period $\leq N$ and thus are
contained in a certain proper analytic subvariety.

\

The result of Hrushovski we are using is as follows:

\begin{theo} (\cite{H}, Corollary 1.2) Let $U$ be an affine variety over
a finite field $\F_q$ and let $S\subset U^2$ be an irreducible subvariety
over $\bar \F_q$. Assume that the two projections of $S$ to $U$ are 
dominant. Denote by $\phi_q$ the Frobenius map. Then for any proper 
subvariety $W$ of $U$, for large enough $m$,
there exists $x\in U(\bar \F_q)$ with $(x, \phi^m_q(x))\in S$ and $x \not\in W$.  
\end{theo}

In particular, let $\bar X$ be any irreducible variety defined over 
a finite field $\F_q$, and 
let
$f:\bar X\dasharrow \bar X$ be a separable dominant rational self-map 
defined over the 
same field.
Let $Y= I\cup R$, where $I$ is the indeterminacy locus and $R$ is the
ramification locus of $f$. Those are subvarieties defined over a finite 
extension of $\F_q$,
and therefore they are periodic under $\phi_q$: for some $k\in \N$ and any
$l\in \N$, 
$\phi_q^l(I_j)=\phi_q^{l+k}(I_j)$ and $\phi_q^l(R_j)=\phi_q^{l+k}(R_j)$ for 
all irreducible components $I_j$ of $I$ and $R_j$ of $R$. Set 
$V=\bar X-Y\cup\phi_q(Y)\cup \dots \cup \phi^k_q(Y)$. Let $U$ be an affine open
subset of $V$ defined over $\F_q$ and let $S$ be the intersection of
 $U\times U$ with the graph
of $f$. Then by Hrushovski's theorem we have the 
following

\begin{coro} \label{overfinite} In the setting as above, there is a point
$x\in \bar X(\bar \F_q)$ such that no iterate $f^i(x)$ is an indeterminacy or
ramification point of $f$, and $x$ is $f$-periodic. Moreover such points
are Zariski-dense in $\bar X$.
\end{coro}

Let now $X$ be a variety defined over a number field $K'$, and let 
$f:X\dasharrow X$ be a dominant rational self-map defined over $K'$. We are going to use
corollary \ref{overfinite} to find, for a suitable finite extension
$K$ of $K'$ and for a suitable prime $\p\subset {\cal O}_K$, a "$\p$-adic 
neighbourhood" in $X$, invariant under some power of $f$, with 
good properties as in \cite{ABR}. The procedure is almost the same as in
\cite{ABR}. Take an affine $U\subset X$ such that $f$ is regular on $U$,
together with a surjective $K'$-morphism $\pi=(x_1, \dots, x_n):U \to \A^n$
(Noether normalization). Write ${\cal O}_U=K'[x_1, \dots, x_n, x_{n+1},
\dots, x_m]/I$, where $I$ contains for instance the minimal polynomials $P_i$ 
of $x_{n+i},\ i>0$ over $K[x_1, \dots, x_n]$ (but probably also something else), so
that $$U=Spec(K'[x_1, \dots, x_n, x_{n+1}, \dots, x_m]/I);$$ we may suppose that
$I$ is given by a system of generators with coefficients from ${\cal O}_{K'}$
and take a model over ${\cal O}_{K'}$: ${\cal U}=
Spec ({\cal O}_{K'} [x_1, \dots, x_n, x_{n+1}, \dots, x_m]/I)$. By abuse of 
notation, we denote the rational map on the model by the same $f$.
Consider $x_{n+1}, \dots x_m, f^*x_1, \dots f^*x_m$ as power series in 
$x_1, \dots, x_n$.
By Lemma 2.1 of \cite{ABR}, their coefficients are $\p'$-integral for almost
all primes $\p'\subset {\cal O}_{K'}$. Choose $\p'$ with this property, and,
moreover, such that the minimal monic polynomials $P_i$ have $\p'$-integral
coefficients and the derivatives $P_i'$ are not identically zero modulo $\p'$, 
and such that $f$ reduced modulo $\p'$ is well-defined and separated. 
Consider the reduction $\bar U$ modulo $\p'$:
$\bar U= Spec(({\cal O}_{K'}/\p')[x_1, \dots, x_n, x_{n+1}, \dots, x_m]/\bar I)$. 
It is equipped with a rational self-map $\bar f$, the reduction of $f$.
By corollary \ref{overfinite}, we can find an $\bar f$-periodic point $x\in \bar U$ over
some finite extension $\F_q(\alpha)$ (where $\F_q={\cal O}_{K'}/\p'$), 
such that no iterate of $x$ by $\bar f$ is in the indeterminacy or ramification,
and such that the values of the derivatives $P_i'$ calculated at $x$ are
non-zero modulo $\p'$. If $\bar U$ is singular, let us moreover choose $x$ in
its smooth locus (however, when our original $U$ is smooth, we can already 
choose the prime $\p'$ in such a way that $\bar U$ is smooth). Let $k$ be the 
period: $\bar f^k(x)=x$.

Let $\beta$ be an algebraic number integral over ${\cal O}_{K'}$ 
such that the reduction of its monic
minimal polynomial over ${\cal O}_{K'}$ modulo $\p'$ gives the minimal 
polynomial 
of $\alpha$,
and let $K=K'(\beta)$. Let $\p$ be some prime of $K$ 
lying over $\p'$.
The point $x$ lifts to a point $y\in U(K_{\p})$ by Hensel's lemma; 
alternatively, we can take as our $(K, \p)$ a slightly larger 
finite extension of 
$(K', \p')$, in order to 
produce a point $y\in U(K)$ which reduces to $x$. 
Define the $\p$-adic neighbourhood
$O_{\p, y}$ of $y$ as follows: 
$$O_{\p, y}=\{t\in U(K_{\p})|x_i(t)\equiv x_i(y) \pmod \p \ for\ 
1\leq i\leq m\}.$$

We may suppose that $\pi(y)=(0,\dots 0)\in \A^n$.
Then, exactly as in \cite{ABR}, we get the following

\begin{prop}\label{nbhd}

(1) The functions $x_1, \dots x_n$ give a bijection
between $O_{{\p}, y}$ and the n-th cartesian power of $\p$.

(2) The set $O_{\p,y}$ contains no indeterminacy and no ramification points 
of $f$.

(3) $f^k(O_{{\p}, y})\subset O_{{\p}, y}$, moreover, $f^k$ is bijective on
$O_{{\p}, y}$.

(4) The $\bar \Q$-points of $X$ are dense in $O_{{\p}, y}$.

\end{prop}

It will be more convenient for us to identify $O_{\p,y}$ with a cartesian
power of ${\cal O}_{\p}$ rather than that of $\p$. So, if our map $f^k$ is given
(say on the completion of the local ring at $y$) by the power series 
$H_i(x_1, \dots x_n)=f^*x_i\in {\cal O}_{\p}[[x_1, \dots x_n]],$ set 
$F(t_1, \dots, t_n)=\frac{1}{r} H(rx_1,\dots rx_n)$ where $r$ is some fixed
 uniformizing element in ${\cal O}_{\p}$ . In this way, we may view
$O_{\p,y}$ as ${\cal O}_{\p}^n$ with coordinates $t_i$, and the map $f^k$ is given by 
the power series $F_i$ on $O_{\p,y}$. Note that, as in \cite{BGT}, the 
$F_i$ have
integral coefficients (by construction the constant terms of $H_i$ are 
divisible by $r$), and moreover the coefficient of $t_1^{k_1}\dots t_n^{k_n}$
in $F_i$ is divisible by $r^{k_1+\dots + k_n-1}$ when $k_1+\dots + k_n\geq 1$
(remark 2.3 of \cite{BGT}). 

Write $F=(F_1,\dots F_n): {\cal O}_{\p}^n\to {\cal O}_{\p}^n$ (so $F$ is a way to think of
the restriction
of $f^k$ to the $\p$-adic neighbourhood) and consider $F$
modulo $\p$. Again as in \cite{BGT}, we have the following 

\begin{prop} \label{idmodp} There is a positive integer $l$ such that 
for every $z\in {\cal O}_{\p}^n$, $F^l(z)\equiv z \pmod \p$.
\end{prop}

{\it Proof:} Modulo $\p$, $F$ is an affine transformation of the linear space
$({\cal O}_{\p}/\p)^n$. Its linear part $L$ is invertible by Proposition 2.4 of
\cite{BGT} since no $f$-iterate of the smooth
point $x$ into which our $\p$-adic neighbourhood reduces modulo $\p$ is
ramification or indeterminacy, and so $\bar f^k$ is unramified at $x$.
Therefore $F$ modulo $\p$ is an automorphism of a finite-dimensional 
affine space
over a finite field, and some power of it is the identity.

\medskip

In Section 3 of \cite{BGT}, the similar situation is considered; the only
difference is that there, the coefficients of all power series are 
in $\Z_p$ rather
than in an extension ${\cal O}_{\p}$. The authors prove:
\begin{theo}\label{main-bgt}(\cite{BGT}, Theorem 3.3)
Let $\phi_1,\dots, \phi_n\in \Z_p[[x_1,\dots, x_n]]$ be convergent power
series on $\Z_p^n$ such that $\phi_i(x) \equiv x_i \pmod p$ and the coefficient 
of 
$x_1^{k_1}\dots x_n^{k_n}$ in the series $\phi_i$ is divisible by 
$p^{k_1+\dots +k_n-1}$ for $k_1+\dots +k_n>1$. Let $(\omega_1,\dots \omega_n)\in
\Z_p^n$. If $p>3$, there exist $p$-adic analytic functions $g_1, \dots g_n\in
\Q_p[[z]]$, convergent on $\Z_p$, such that $g_i(\Z_p)\subset \Z_p$,
$g_i(0)=\omega_i$ and $g_i(z+1)=\phi_i(g_1(z),\dots , g_n(z))$.
\end{theo}
They construct $g_i(z)$ "by approximation", as a Mahler series 
$$g_i(z)=\omega_i+\sum_{k=1}^\infty b_{ik}\binom{z}{k},$$ and it turns out,
by their construction, that $b_{ik}$ are of the 
form
$\sum_{j=(k+1)/2}^\infty p^jc_{ijk}$, where $c_{ijk}\in \Z_p$. Since 
$|b_{ik}|_p\to 0$ when $k\to \infty$, these Mahler series define continuous
functions on $\Z_p$ with values in $\Z_p$. To show that these functions are
in fact analytic on $\Z_p$, one needs to check (\cite{R}, Theorem 4.7 of ChapterVI) that
$|b_{ik}|_p/|k!|_p\to 0$ when $k\to \infty$, and this is true for $p>3$ since
$|b_{ik}|_p\leq p^{-(k+1)/2}$ and $1/|k!|_p<p^{k/(p-1)}$.

\medskip

In the situation when $\Phi=(\phi_1, \dots , \phi_n)$ is an analytic map 
from ${\cal O}_{\p}^n$ to itself and so $\Z_p$, not as the domain of definition of $g_i$ 
but as the domain where the $g_i$ take their values, is replaced by 
the extension ${\cal O}_{\p}$,
their argument goes through almost verbatim, provided that one replaces $p$ with a uniformizing 
element $r$ of ${\cal O}_{\p}$ where appropriate: indeed the theory of Mahler 
series applies
to ${\cal O}_{\p}$-valued functions on $\Z_p$ as well (\cite{R}, 
chapters IV.2.3,VI.4.7). The only 
exception is the last step concerning the analyticity of the ${\cal O}_{\p}$-
valued map $g_i$ on $\Z_p$: instead of being analytic on the whole of
$\Z_p$, it is going to be analytic on a certain neighbourhood $p^l\Z_p$.
This is because in the expression $b_{ik}=\sum_{j=(k+1)/2}^\infty p^jc_{ijk},$
we have to replace $p$ by the uniformizing element $r$, which can be of 
smaller $p$-adic order $1/e$; so that $|b_{ik}|_p\leq p^{-(k+1)/2e}$ and 
we need the condition $p>2(e+1)$ in order to guarantee $|b_{ik}|_p/|k!|_p\to 0$.
But one always has convergence on some $p^l\Z_p$. Indeed, one has the
following elementary lemma:

\begin{lemm}\label{convergence} A Mahler series $\sum_{k=1}^\infty b_{k}\binom{z}{k}$
defines an analytic function on $p^l\Z_p$ as soon as
$$|p^{k\frac{p^l-1}{(p-1)p^l}}b_{k}/k!|_p\to 0.$$
\end{lemm}

{\it Proof:} Write each term $\binom{p^lz}{k}$ as a product of a constant
$C_k$ and a polynomial with integral coefficients. As soon as $b_kC_k\to 0$, our 
Mahler series will be a convergent power series on $p^l\Z_p$ and thus analytic
on $p^l\Z_p$. A standard computation
shows that one can take
$$C_k=\frac{1}{k!}p^lp^{[\frac{k-1}{p}]+[\frac{1}{p}[\frac{k-1}{p}]]+\dots+
[\frac{1}{p}[\dots[\frac{k-1}{p}]\dots]]},$$
where $[x]$ denotes the integral part of $x$ and in the last summand, the division
by $p$ is made $l$ times. The condition $b_kC_k\to 0$ is equivalent to the one
in the statement of our lemma: indeed, since we are only interested in the terms which 
grow together with $k$, we can forget about the integral parts and use the
formula for the sum of a geometric progression. 

\medskip

Clearly, the condition in the above lemma is satisfied by the Mahler series 
$\sum_{k=1}^\infty b_{ik}\binom{z}{k}$ as soon as $l$ is not too small with respect
to $e$.

So the theorem of \cite{BGT} becomes

\begin{theo}\label{param}Let $\phi_1,\dots, \phi_n\in {\cal O}_{\p}[[x_1,\dots, x_n]]$ be convergent power
series such that $\phi_i(x)\equiv x_i \pmod \p$ and the coefficient of 
$x_1^{k_1}\dots x_n^{k_n}$ in the series $\phi_i$ is divisible by 
$r^{k_1+\dots +k_n-1}$ for $k_1+\dots +k_n>1$ for a uniformizing element $r$.  
Let $(\omega_1,\dots \omega_n)\in
{\cal O}_{\p}^n$. Then there exist functions $g_1, \dots g_n$, continious on
$\Z_p$ and
analytic on $p^l\Z_p$ for a certain positive integer $l$, 
such that $g_i(\Z_p)\subset {\cal O}_{\p}$,
$g_i(0)=\omega_i$ and $g_i(z+1)=\phi_i(g_1(z),\dots , g_n(z))$.
\end{theo}

Here is an immediate corollary of Theorem \ref{param}:

\begin{coro}\label{boundedness} 
Let $X$, $f:X\dasharrow X$ be a variety and a dominant rational self-map defined over
a number field. Let $O_{\p, y}$ be an $f^k$-invariant $\p$-adic neighbourhood 
constructed in proposition \ref{nbhd} (formed by $K_{\p}$-points for a suitable extension
$K$ of our number field and a suitable prime $\p$).  
Then there exists a positive integer $N$ such that any preperiodic point 
in this neighbourhood
is periodic of bounded period $\leq N$.
\end{coro}

{\it Proof:} Indeed, the neighbourhood is $F=f^k$-invariant, and some further
power $F^l$ of $f$ satisfies the conditions for $\phi$ in theorem \ref{param}.
Let $\omega$ be a preperiodic point. From theorem \ref{param}, we deduce that its orbit
under a still larger power $\psi=\phi^s=f^N$, depending only on the 
neighbourhood
itself (in particular, on the ramification index of $K_{\p}$ over $\Q_p$),
is encoded by an analytic map $g:p^l\Z_p\to O_{\p, y}$, satisfying
$g(p^li)=\psi^i(\omega)$ for $i\in \N$. But an analytic function taking some
of its values infinitely many times on $p^l\Z$ must be constant, so in 
fact $\psi(\omega)=\omega$
and so all points with a finite orbit in the neighbourhood are $N$-periodic. 

\medskip

From this, we immediately draw the conclusion announced in the beginning:

\begin{coro}\label{nonpreper} Unless if $f$ is of finite order, 
there exist points in $X(\bar \Q)$ which are
not preperiodic with respect to $f$.
\end{coro}

{\it Proof:} Indeed, points of bounded period must be contained in a 
proper analytic subvariety of $O_{\p, y}$, whereas algebraic points are
dense in $O_{\p, y}$.

\medskip

{\bf Acknowledgements:} This note has been written during the author's
stay at the Institute for Advanced Study. The author has been supported
by the NSF grant DMS-0635607 and by the Minerva Research Foundation. 
I thank the IAS
for its hospitality and excellent working conditions, and NSF and 
Minerva for the support which made my stay possible.

I am grateful to Marat Rovinsky, Antoine Chambert-Loir, P\"ar Kurlberg,
Dragos Ghioca and Tom Tucker for helpful discussions.

\end{document}